\documentclass[a4paper,twoside]{article}

\makeatletter
\renewcommand{\@seccntformat}[1]{\csname the#1\endcsname}
\renewcommand\normalsize{%
   \@setfontsize\normalsize\@xpt\@xiipt
   \abovedisplayskip 6\p@ \@plus1\p@
   \abovedisplayshortskip 6\p@ \@plus1\p@
   \belowdisplayshortskip 6\p@ \@plus1\p@
   \belowdisplayskip \abovedisplayskip
   \let\@listi\@listI}
\normalsize
\renewcommand\ps@headings{%
  \let\@oddfoot\@empty\let\@evenfoot\@empty
  \let\@oddhead\@empty\let\@evenhead\@empty
  \let\@mkboth\markboth
  \def\sectionmark##1{%
    \markboth {\uppercase{\ifnum \c@secnumdepth >\z@
        \thesection.\relax\fi
        ##1}}{}}%
  \def\subsectionmark##1{%
    \markright {\ifnum \c@secnumdepth >\@ne
        \thesubsection\relax \fi
        ##1}}}
\renewcommand\ps@myheadings{%
    \let\@oddfoot\@empty\let\@evenfoot\@empty
    \def\@evenhead{\thepage\hfil\small\leftmark\hfil}%
    \def\@oddhead{\hfil{\small\rightmark}\hfil\thepage}%
    \let\@mkboth\@gobbletwo
    \let\sectionmark\@gobble
    \let\subsectionmark\@gobble
    }
\def\myfnsymbol#1{\expandafter\@myfnsymbol\csname c@#1\endcsname}
\def\@myfnsymbol#1{\ensuremath{\ifcase#1\or *\or \ddagger\or **\or
   \mathsection\or \mathparagraph\or \|\or **\or \dagger\dagger
   \or \ddagger\ddagger \else\@ctrerr\fi}}
\renewcommand\maketitle{\par
  \begingroup
    \renewcommand\thefootnote{\myfnsymbol{footnote}}%
    \def\@makefnmark{\hbox to\z@{$\m@th^{\@thefnmark}$\hss}}%
    \long\def\@makefntext##1{\parindent 1em\noindent
            \hbox to1.8em{\hss$\m@th^{\@thefnmark}$}##1}%
    \if@twocolumn
      \ifnum \col@number=\@ne
        \@maketitle
      \else
        \twocolumn[\@maketitle]%
      \fi
    \else
      \newpage
      \global\@topnum\z@   % Prevents figures from going at top of page.
      \@maketitle
    \fi
    \thispagestyle{headings}
    \@thanks
  \endgroup
  \setcounter{footnote}{0}%
  \let\thanks\relax
  \let\maketitle\relax\let\@maketitle\relax
  \gdef\@thanks{}\gdef\@author{}\gdef\@title{}}
\renewcommand\@maketitle{%
  \newpage
  \begin{center}%
    {\LARGE\bfseries \@title \par}%
    \vskip 24\p@%
        {\Large\itshape
      \lineskip .5em%
      \begin{tabular}[t]{c}%
        \@author
      \end{tabular}\par}%
  \end{center}%
  \par
  \vskip 60\p@}
\def\@sect#1#2#3#4#5#6[#7]#8{\ifnum #2>\c@secnumdepth
     \let\@svsec\@empty\else
     \refstepcounter{#1}%
     \let\@@protect\protect
     \def\protect{\noexpand\protect\noexpand}%
     \edef\@svsec{\@seccntformat{#1}}%
     \let\protect\@@protect\fi
     \@tempskipa #5\relax
      \ifdim \@tempskipa>\z@
      {#6\relax{\interlinepenalty \@M \@svsec #8\par}}%
               \csname #1mark\endcsname{#7}\addcontentsline
         {toc}{#1}{\ifnum #2>\c@secnumdepth \else
                      \protect\numberline{\csname the#1\endcsname}\fi
                    #7}\else
        \def\@svsechd{#6%
\hskip #3\relax  %% \relax added 2 May 90
                   \@svsec #8\csname #1mark\endcsname
                      {#7}\addcontentsline
                           {toc}{#1}{\ifnum #2>\c@secnumdepth \else
                           \protect\numberline{\csname the#1\endcsname}%
                                     \fi
                       #7}}\fi
     \@xsect{#5}}
\def\@ssect#1#2#3#4#5{\@tempskipa #3\relax
   \ifdim \@tempskipa>\z@
     {#4%
        {\interlinepenalty \@M #5\par}%
     }
   \else \def\@svsechd{#4%
\hskip #1
    \relax
                        #5}\fi
    \@xsect{#3}}

\renewcommand\section{\@startsection {section}{1}{\z@}%
                                  {-36\p@ \@plus -1\p@ \@minus -4\p@}%
                                  {12\p@ \@plus 1\p@}%
{\reset@font\Large\bfseries\centering}}
\renewcommand\subsection{\@startsection{subsection}{2}{\z@}%
                                     {-24\p@ \@plus -1\p@ \@minus -4\p@}%
                                     {8\p@ \@plus 1\p@}%12
                                     {\reset@font\large\bfseries\centering}}

\renewcommand\subsubsection{\@startsection{subsubsection}{3}{\z@}%
                                     {3.25ex \@plus 1ex \@minus .2ex}%
                                     {-1.5ex \@plus .2ex}%
                                     {\reset@font\normalsize\bfseries}}

\renewenvironment{abstract}{\noindent\small{\bfseries\abstractname.}%
        }%
        {\vskip 11\p@}
\newenvironment{classification}{\noindent\small 2000 Mathematics 
Subject Classification:}{\vskip 12\p@}
\def\@lbibitem[#1]#2{\item[\hfill\@biblabel{#1}]\if@filesw
      {\let\protect\noexpand
       \immediate
       \write\@auxout{\string\bibcite{#2}{#1}}}\fi\ignorespaces}

\renewenvironment{thebibliography}[1]
     {\section*{\reset@font\fontsize{11.6}{13.6pt}\bfseries\refname
        \@mkboth{\uppercase{\refname}}{\uppercase{\refname}}}%
      \list{\@biblabel{\arabic{enumiv}}}%
           {\settowidth\labelwidth{\@biblabel{#1}}%
            \setlength{\leftmargin}{\labelwidth}
            \labelsep2mm
            \itemsep4pt
        \topsep\z@
            \parsep\z@
            \advance\leftmargin\labelsep
                        \usecounter{enumiv}%
            \let\p@enumiv\@empty
            }%
            \sloppy\clubpenalty4000\widowpenalty4000%
      \sfcode`\.=\@m
      \small}
     {\def\@noitemerr
       {\@latex@warning{Empty `thebibliography' environment}}%
      \endlist}

\def\@seccntformat#1{\csname the#1\endcsname.\hskip 0.25em}
\newskip\aline \newskip\halfaline
\def\skipaline{\vskip\aline}

\def\qedbox{$\rlap{$\sqcap$}\sqcup$}

\def\Proof{\ifdim\lastskip<\aline\removelastskip\skipaline\fi
\noindent\it Proof. \rm}

\RequirePackage{amsmath}
\RequirePackage{amssymb}
\RequirePackage{ifthen}
\RequirePackage{xypic}

\def\FMithmInfo{1995/11/23 v2.2c Theorem extension package (FMi)}
\@ifundefined{theorem@style}{}{\endinput}

\gdef\theoremstyle#1{%
   \@ifundefined{th@#1}{\@warning
          {Unknown theoremstyle `#1'. Using `plain'}%
          \theorem@style{plain}}%
      {\theorem@style{#1}}%
      \begingroup
        \csname th@\the\theorem@style \endcsname
      \endgroup}
\global\let\@begintheorem\relax
\global\let\@opargbegintheorem\relax
\newtoks\theorem@style
\global\theorem@style{plain}
\gdef\theorembodyfont#1{%
   \def\@tempa{#1}%
   \ifx\@tempa\@empty
    \theorem@bodyfont{}%
   \else
    \theorem@bodyfont{\reset@font#1}%
   \fi
   }
\newtoks\theorem@bodyfont
\global\theorem@bodyfont{}
\gdef\theoremheaderfont#1{\gdef\theorem@headerfont{#1}%
       \gdef\theoremheaderfont##1{%
        \typeout{\string\theoremheaderfont\space should be used
                 only once.}}}
\ifx\upshape\undefined
\gdef\theorem@headerfont{\bfseries}
\else \gdef\theorem@headerfont{\normalfont\bfseries}\fi
\gdef\th@plain{\@input@{thp.sty}}
\gdef\th@break{\@input@{thb.sty}}
\gdef\th@marginbreak{\@input@{thmb.sty}}
\gdef\th@changebreak{\@input@{thcb.sty}}
\gdef\th@change{\@input@{thc.sty}}
\gdef\th@changep{\@input@{thcp.sty}}% only new line, compared with theorem.sty
\gdef\th@changepd{\@input@{thcpd.sty}}% only new line, compared with theorem.sty
\gdef\th@margin{\@input@{thm.sty}}
\gdef\@xnthm#1#2[#3]{\expandafter\@ifdefinable\csname #1\endcsname
   {%
    \@definecounter{#1}\@newctr{#1}[#3]%
    \expandafter\xdef\csname the#1\endcsname
      {\expandafter \noexpand \csname the#3\endcsname
       \@thmcountersep \@thmcounter{#1}}%
    \def\@tempa{\global\@namedef{#1}}%
    \expandafter \@tempa \expandafter{%
      \csname th@\the \theorem@style
            \expandafter \endcsname \the \theorem@bodyfont
     \@thm{#1}{#2}}%
    \global \expandafter \let \csname end#1\endcsname \@endtheorem
   }}
\gdef\@ynthm#1#2{\expandafter\@ifdefinable\csname #1\endcsname
   {\@definecounter{#1}%
    \expandafter\xdef\csname the#1\endcsname{\@thmcounter{#1}}%
    \def\@tempa{\global\@namedef{#1}}\expandafter \@tempa
     \expandafter{\csname th@\the \theorem@style \expandafter
     \endcsname \the\theorem@bodyfont \@thm{#1}{#2}}%
    \global \expandafter \let \csname end#1\endcsname \@endtheorem}}
\gdef\@othm#1[#2]#3{%
  \expandafter\ifx\csname c@#2\endcsname\relax
   \@nocounterr{#2}%
  \else
   \expandafter\@ifdefinable\csname #1\endcsname
   {\expandafter \xdef \csname the#1\endcsname
     {\expandafter \noexpand \csname the#2\endcsname}%
    \def\@tempa{\global\@namedef{#1}}\expandafter \@tempa
     \expandafter{\csname th@\the \theorem@style \expandafter
     \endcsname \the\theorem@bodyfont \@thm{#2}{#3}}%
    \global \expandafter \let \csname end#1\endcsname \@endtheorem}%
  \fi}
\gdef\@thm#1#2{\refstepcounter{#1}%
   \trivlist
   \@topsep \theorempreskipamount               % used by first \item
   \@topsepadd \theorempostskipamount           % used by \@endparenv
   \@ifnextchar [%
   {\@ythm{#1}{#2}}%
   {\@begintheorem{#2}{\csname the#1\endcsname}\ignorespaces}}
\global\let\@xthm\relax
\newskip\theorempreskipamount
\newskip\theorempostskipamount
\global
\global
\global\let\@endtheorem=\endtrivlist
\@onlypreamble\@xnthm
\@onlypreamble\@ynthm
\@onlypreamble\@othm
\@onlypreamble\newtheorem
\@onlypreamble\theoremstyle
\@onlypreamble\theorembodyfont
\@onlypreamble\theoremheaderfont
\theoremstyle{plain}

\begingroup \makeatletter
\@ifundefined{theorem@style}{\input{theorem.sty}}{}
             [\FMithmInfo]
\gdef\th@plain{\normalfont\itshape
  \def\@begintheorem##1##2{%
        \item[\hskip\labelsep \theorem@headerfont ##1\ ##2.]}%
\def\@opargbegintheorem##1##2##3{%
   \item[\hskip\labelsep \theorem@headerfont ##1\ ##2\ \rm(##3).]}}
\endgroup

\def\@addpunct#1{\ifnum\spacefactor>\@m \else#1\fi}

\newenvironment{proof}[1][\proofname]{\par
  \normalfont
  \topsep6\p@\@plus6\p@ \trivlist
  \item[\hskip\labelsep\itshape
    #1\@addpunct{.}]\ignorespaces
}{\proof@ending\endtrivlist\par}

\newcommand{\proofname}{Proof}
\newcommand{\proof@ending}{\hfill\qedbox} \newcommand{\BoxedProofs}{\renewcommand{\proof@ending}{\hfill\ensuremath{\Box}}}
\newcommand{\NonBoxedProofs}{\renewcommand{\proof@ending}{\hfill\qedbox}}

\newskip\aline \newskip\halfaline
\def\skipaline{\vskip\aline}

\newcommand{\@lstlabel}{}
\newcommand{\lstlabel}[1]{\renewcommand{\@lstlabel}{#1}}
\newcommand{\@lsttemplate}{}
\newcommand{\lsttemplate}[1]{\renewcommand{\@lsttemplate}{#1}}
{\renewcommand{\theenumi}{\@lstlabel}%
  \begin{list}%
    {\@lstlabel}%
    {\usecounter{enumi}%
      \settowidth{\labelwidth}{\@lsttemplate}%
      \setlength{\leftmargin}{\labelwidth}%
      \setlength{\labelsep}{0pt}%
      \setlength{\topsep}{\medskipamount}%
      \setlength{\parsep}{0pt}%
      \setlength{\itemsep}{0pt}%
      \setlength{\itemindent}{0pt}%
      \setlength{\listparindent}{\parindent}}}
  {\end{list}}
\newcommand{\deflststyle}[2]{\expandafter\newcommand\expandafter{\csname @@#1@@\endcsname}{#2}}
\newcommand{\lststyle}[1]{\csname @@#1@@\endcsname}
\deflststyle{}{
  \lsttemplate{}
  \lstlabel{}}
\deflststyle{ }{
  \lsttemplate{\hspace{\parindent}}
  \lstlabel{}}
\deflststyle{--}{
  \lsttemplate{\hspace{\parindent}}
  \lstlabel{\,--\hfill}}
\deflststyle{a.}{
  \lsttemplate{\textnormal{\,b.\ \ }}
  \lstlabel{\textnormal{\,\alph{enumi}.\hfill}}}
\deflststyle{a)}{
  \lsttemplate{\textnormal{\,b). }}
  \lstlabel{\textnormal{\,\alph{enumi})\hfill}}}
\deflststyle{(a)}{
  \lsttemplate{\textnormal{(b)\ \ }}
  \lstlabel{\textnormal{(\alph{enumi})\hfill\ \ }}}
\deflststyle{1.}{
  \lsttemplate{\textnormal{\,0.\ \ }}
  \lstlabel{\textnormal{\hfill\arabic{enumi}.\ \ }}}
\deflststyle{1)}{
  \lsttemplate{\textnormal{\,0)\ \ }}
  \lstlabel{\textnormal{\hfill\arabic{enumi})\ \ }}}
\deflststyle{(1)}{
  \lsttemplate{\textnormal{(0)\ \ }}
  \lstlabel{\textnormal{\hfill(\arabic{enumi})\ \ }}}
\deflststyle{i.}{
  \lsttemplate{\textnormal{\,iii.\ \ }}
  \lstlabel{\textnormal{\hfill\roman{enumi}.\ \ }}}
\deflststyle{i)}{
  \lsttemplate{\textnormal{\,iii)\ \ }}
  \lstlabel{\textnormal{\hfill\roman{enumi})\ \ }}}
\deflststyle{(i)}{
  \lsttemplate{\textnormal{(iii)\ \ }}
  \lstlabel{\textnormal{\hfill(\roman{enumi})\ \ }}}
\lststyle{(a)} % default style

\newenvironment{block}{\begin{trivlist}\item{}}{\end{trivlist}}

\reversemarginpar

\let\s@bsection=\subsection
\newcommand{\c@mda}[2][]{\s@bsection[#1]{#2.}}
\newcommand{\c@mdb}[1]{\s@bsection*{#1.}}
\renewcommand{\subsection}{\secdef\c@mda\c@mdb}

\let\oldcite=\cite
\renewcommand{\cite}[2][no@ption]
{\ifthenelse{\equal{#1}{no@ption}}{\textnormal{\oldcite{#2}}}{\oldcite[#1]{#2}}}
\let\oldref=\ref
\renewcommand{\ref}[1]{\textnormal{\oldref{#1}}}
\let\oldpageref=\ref
\renewcommand{\pageref}[1]{\textnormal{\oldpageref{#1}}}

\let\@contact=\empty
\newcommand{\contact}[2][]{
  \expandafter\gdef\expandafter\@contact\expandafter{%
    \@contact
  \small
  \begin{block}#2\\[.5em]
  \ifthenelse{\equal{#1}{}}{}{Email: #1}\end{block}
  }}
\newcommand{\makelastpage}{\medskip\@contact}

\DeclareMathOperator{\id}{id}

\newcommand{\BDC}{\mathbf{D}^{\mathrm{b}}}

\newcommand{\ilim}[1][]{\mathop{\varinjlim}\limits_{#1}}

\renewcommand{\to}[1][]{\xrightarrow[#1]{}}
\newcommand{\from}[1][]{\xleftarrow[#1]{}}
\newcommand{\isoto}[1][]{\xrightarrow[#1]{\sim}}

\newcommand{\Endo}[1][]{\mathrm{End}_{\raise1.5ex\hbox to.1em{}#1}}

\newcommand{\Hom}[1][]{\mathrm{Hom}_{\raise1.5ex\hbox to.1em{}#1}}

\newcommand{\RHom}[1][]{\mathrm{RHom}_{\raise1.5ex\hbox to.1em{}#1}}

\newcommand{\Ext}[2][]{\mathrm{Ext}_{\raise1.5ex\hbox to.1em{}#1}^{#2}}

\newcommand{\THom}[1][]{\mathrm{THom}_{\raise1.5ex\hbox to.1em{}#1}}

\newcommand{\Mod}{\mathrm{Mod}}

\newcommand{\Tens}[1][]{\mathbin{\otimes_{\raise1.5ex\hbox to-.1em{}#1}}}

\newcommand{\LTens}[1][]{\mathbin{\otimes_{\raise1.5ex\hbox to-.1em{}#1}^{L}}}

\newcommand{\Tor}[2][]{\mathrm{Tor}^{\raise1.5ex\hbox to.1em{}#1}_{#2}}

\def\shf{\mathcal{F}}

\def\shi{\mathcal{I}}

\def\shl{\mathcal{L}}
\def\shm{\mathcal{M}}
\def\shn{\mathcal{N}}

\def\shp{\mathcal{P}}
\def\shq{\mathcal{Q}}

\newcommand{\sect}{\varGamma}
\newcommand{\rsect}{\mathrm{R}\varGamma}

\renewcommand{\hom}[1][]{{\mathcal{H}om}_{\raise1.5ex\hbox to.1em{}#1}}

\newcommand{\rhom}[1][]{{R\mathcal{H}om}_{\raise1.5ex\hbox to.1em{}#1}}

\newcommand{\ext}[2][]{{\mathcal{E}xt}_{\raise1.5ex\hbox to.1em{}#1}^{#2}}

\newcommand{\thom}[1][]{{T\mathcal{H}om}_{\raise1.5ex\hbox to.1em{}#1}}

\newcommand{\tens}[1][]{\mathbin{\otimes_{\raise1.5ex\hbox to-.1em{}#1}}}

\newcommand{\ltens}[1][]{\mathbin{\otimes_{\raise1.5ex\hbox to-.1em{}#1}^{L}}}

\newcommand{\tor}[2][]{{\mathcal{T}or}^{\raise1.5ex\hbox to.1em{}#1}_{#2}}

\newcommand\etens{\mathbin{\boxtimes}}

\newcommand{\roim}[1]{{R#1}_*}

\newcommand{\opb}[1]{#1^{-1}}

\newcommand{\GHom}[1][]{\mathrm{GHom}_{\raise1.5ex\hbox to.1em{}#1}}

\newcommand{\GExt}[2][]{\mathrm{GExt}_{\raise1.5ex\hbox to.1em{}#1}^{#2}}

\newcommand{\FHom}[1][]{\mathrm{FHom}_{\raise1.5ex\hbox to.1em{}#1}}

\newcommand{\ghom}[1][]{{\mathcal{GH}om}_{\raise1.5ex\hbox to.1em{}#1}}

\newcommand{\gext}[2][]{{\mathcal{GE}xt}_{\raise1.5ex\hbox to.1em{}#1}^{#2}}

\newcommand{\fhom}[1][]{{\mathcal{FH}om}_{\raise1.5ex\hbox to.1em{}#1}}

\newcommand{\tenstop}[1][]{\mathbin{\hat{\otimes}_{\raise1.5ex\hbox to-.1em{}#1}}}

\newcommand{\homtop}[1][]{\mathcal{L}_{\raise1.5ex\hbox to.1em{}#1}}

\newcommand{\Homtop}[1][]{\mathrm{L}_{\raise1.5ex\hbox to.1em{}#1}}

\newcommand{\D}{\mathcal{D}}

\renewcommand{\O}{\mathcal{O}}

\def\absdoim#1{\underline{#1}_*}
\def\reldoim[#1]#2{\underline{#2}_{|{#1}*}}
\def\doim{\@ifnextchar [{\reldoim}{\absdoim}}

\def\absdeim#1{\underline{#1}_*}
\def\reldeim[#1]#2{\underline{#2}_{|{#1}*}}
\def\deim{\@ifnextchar [{\reldeim}{\absdeim}}

\def\absdopb#1{\underline{#1}^{-1}}
\def\reldopb[#1]#2{\underline{#2}_{|{#1}}^{-1}}
\def\dopb{\@ifnextchar [{\reldopb}{\absdopb}}

\def\absboim#1{\underline{\underline{#1}}_*}
\def\relboim[#1]#2{\underline{\underline{#2}}_{|{#1}*}}
\def\boim{\@ifnextchar [{\relboim}{\absboim}}

\def\absbeim#1{\underline{\underline{#1}}_*}
\def\relbeim[#1]#2{\underline{\underline{#2}}_{|{#1}*}}
\def\beim{\@ifnextchar [{\relbeim}{\absbeim}}

\def\absbopb#1{\underline{\underline{#1}}^*}
\def\relbopb[#1]#2{\underline{\underline{#2}}_{|{#1}}^*}
\def\bopb{\@ifnextchar [{\relbopb}{\absbopb}}

\newtheorem{theorem}{Theorem}[section]

\theorembodyfont{\rmfamily}

\newtheorem{remark}[theorem]{Remark}

\numberwithin{equation}{section}

\renewcommand{\to}[1][]{\xrightarrow{#1}}
\renewcommand{\from}[1][]{\xleftarrow{#1}}

\renewcommand{\BDC}{\mathsf{D}^{\mathrm{b}}}
\newcommand{\shift}[1]{{[#1]}}

\newcommand{\qcoh}{{\mathrm{qc}}}

\renewcommand{\doim}[1]{#1_+}
\renewcommand{\dopb}[1]{#1^*}
\newcommand{\dtens}{\mathbin{\otimes}}

\newcommand{\OX}{\O_X}
\newcommand{\DX}{\D_X}

\newcommand{\tr}[1]{{}^t \! {#1}}
\newcommand{\ttr}[1]{(^t \! {#1})}

\newcommand{\AffX}{\mathbb{A}^1_X}
\newcommand{\DAffX}{\D_{\AffX}}

\newcommand{\OY}{\O_Y}
\newcommand{\DY}{\D_Y}

\newcommand{\VV}{V}
\newcommand{\piV}{\pi}

\newcommand{\OV}{\O_{\VV}}
\newcommand{\DV}{\D_{\VV}}

\renewcommand{\AA}{\check \VV}
\newcommand{\piA}{\check\pi}
\newcommand{\iA}{\check\iota}
\newcommand{\OA}{\O_{\AA}}
\newcommand{\DA}{\D_{\AA}}

\newcommand{\VVAA}{\VV\times_X\AA}

\newcommand{\WW}{W}

\newcommand{\DW}{\D_{\WW}}

\newcommand{\MM}{\check \WW}

\newcommand{\WWMM}{\WW\times_X\MM}

\newcommand{\VVMM}{\VV\times_X\MM}

\newcommand{\dfourier}[1][]{\shf_{#1}}

\makeatother

\contact[baldassa@math.unipd.it, dagnolo@math.unipd.it]{Dipartimento
di Matematica Pura ed Applicata; Universit\`a di Padova; via G.
Belzoni, 7; 35131 Padova; Italy}

\begin{document}

\author{Francesco Baldassarri \and Andrea D'Agnolo} 

\title{On Dwork cohomology and algebraic $\D$-modules}

\maketitle

\begin{abstract}
After works by Katz, Monsky, and Adolphson-Sperber, a comparison theorem between relative de Rham cohomology and Dwork cohomology is established in a paper by Dimca-Maaref-Sabbah-Saito in the framework of algebraic $\D$-modules. 
We propose here an alternative proof of this result. The use of Fourier transform techniques makes our approach more functorial. 
\end{abstract}


%\tableofcontents

\section{Review of algebraic $\D$-modules}

For the reader's convenience, we recall here the notions and results
from the theory of algebraic $\D$-modules that we need. Our references were \cite{Bo,B,K0,K1}.

\subsection{Basic operations}\label{se:algDmod}

Let $X$ be a smooth algebraic variety over a field of
characteristic zero, and let $\OX$ and $\DX$ be its structure sheaf
and the sheaf of differential operators, respectively.  Let $\Mod(\DX)$
be the abelian category of left $\DX$-modules, $\BDC(\DX)$ its bounded
derived category, and $\BDC_\qcoh(\DX)$ the full triangulated
subcategory of $\BDC(\DX)$ whose objects have quasi-coherent
cohomologies.

Let $f\colon X\to Y$ be a morphism of smooth algebraic varieties, and
denote by $\D_{X\rightarrow Y}$ and $\D_{Y\leftarrow X}$ the transfer bimodules.
We use the following notation for the operations of tensor product,
inverse image, and direct image for $\D$-modules\footnote{About the tensor product, note that $\shm \ltens[\OX] \shm' \simeq (\shm \tens[\OX] \DX) \ltens[\DX] \shm' \simeq \shm \ltens[\DX] (\DX \tens[\OX] \shm')$, where $\shm
\tens[\OX] \DX$ (resp. $\DX \tens[\OX] \shm'$) is given the natural
stucture of left-right (resp. left-left) $\DX$-bimodule, and
$\ltens[\DX]$ always uses up the ``trivial'' $\DX$-module structure.}
\begin{align*}
\dtens\ &\colon \BDC_\qcoh(\DX) \times \BDC_\qcoh(\DX) \to \BDC_\qcoh(\DX), & (\shm ,
\shm') &\mapsto \shm \ltens[\OX] \shm', \\
\dopb f &\colon \BDC_\qcoh(\DY) \to \BDC_\qcoh(\DX), & \shn &\mapsto \D_{X\rightarrow Y}
\ltens[\opb f \DY] \opb f \shn, \\
\doim f &\colon \BDC_\qcoh(\DX) \to \BDC_\qcoh(\DY), & \shm &\mapsto \roim
f(\D_{Y\leftarrow X} \ltens[\DX] \shm).
\end{align*}

If $f\colon X\to Y$ and $g\colon Y\to Z$ are morphisms of smooth algebraic varieties, then there are natural functorial isomorphisms
\begin{eqnarray}
\label{eq:fcircgopb}
\dopb f \dopb g &\simeq& \dopb{(g\circ f)}, \\
\label{eq:fcircgoim}
\doim g \doim f &\simeq& \doim{(g\circ f)}.
\end{eqnarray}
For $\shn,\shn'\in\BDC_\qcoh(\DY)$, there
is a natural isomorphism
\begin{equation}
\label{eq:fopbtens}
\dopb f (\shn \dtens \shn') \simeq \dopb f \shn \dtens \dopb f
\shn'.
\end{equation}
For $\shm\in\BDC_\qcoh(\DX)$ and $\shn\in\BDC_\qcoh(\DY)$, there
is a projection formula\footnote{\label{footnotenum}In the appendix we recall the proofs of base change and projection formulae.}
\begin{equation}
\label{eq:projection}
\doim f(\shm \dtens \dopb f \shn) \simeq \doim f \shm \dtens \shn.
\end{equation}
Consider a Cartesian square of smooth algebraic varieties
\begin{equation}
\label{eq:cartsq}
\xymatrix{ X' \ar[r]^{h'} \ar[d]_{f' }& X \ar[d]_{f}
\ar@{}[dl]|-\square \\ Y' \ar[r]^{h} & Y. }
\end{equation}
For $\shm\in\BDC_\qcoh(\DX)$, there is a base change formula$^\text{\ref{footnotenum}}$
\begin{equation}
\label{eq:basechange}
\doim{f'} \dopb{h'{}} \shm \shift{d_{X'}-d_X}
\simeq 
\dopb h \doim f \shm \shift{d_{Y'}-d_Y},
\end{equation}
where $d_X$ denotes the dimension of $X$.

\subsection{Relative cohomology}

Let $S$ be a closed subscheme of $X$, and denote by
$\shi_{S}\subset\OX$ the corresponding ideal of $\OX$. 
For $\shf\in\Mod(\OX)$ one sets\footnote{In other words, for any open subset $V\subset X$, $\sect_{[S]}
\shf(V) = \{ s\in\shf(V) \colon (\shi_{S}|_V)^{m} s = 0, \ m\gg 0 \}$.  Recall that
if $\shf$ is quasi-coherent, Hilbert's Nullstellensatz implies that
$\sect_{[S]} \shf \simeq \sect_{S} \shf$, the subsheaf of $\shf$ whose
sections are supported in $S$.}
$$
\sect_{[S]}(\shf) = \ilim[m] \hom[\OX] (\OX/\shi_{S}^{m} , \shf).
$$
We point out that $\sect_{[S]}(\shf) = \sect_{[S^{\text{red}}]}(\shf)$.
If $\shm \in \Mod(\DX)$ one checks that $\sect_{[S]} \shm$ has a
natural left $\DX$-module structure, and one considers the right derived
functor
$$
\rsect_{[S]} \colon \BDC_\qcoh(\DX) \to \BDC_\qcoh(\DX).
$$

Let $i\colon X\setminus S \to X$ be the open embedding, and $\shm\in\BDC_\qcoh(\DX)$. There is
a distinguished triangle in $\BDC_\qcoh(\DX)$
\begin{equation}
\label{eq:rsectDT}
\rsect_{[S]}\shm \to \shm \to \doim i \dopb i \shm \to[+1].
\end{equation}
For $S,S'\subset X$ possibly singular closed subvarieties, and $\shm\in\BDC_\qcoh(\DX)$, one has
\begin{eqnarray}
\label{eq:sectM}
\rsect_{[S]} \shm &\simeq& \shm \dtens \rsect_{[S]} \OX, \\
\label{eq:sectSS'}
\rsect_{[S]} \rsect_{[S']} \shm &\simeq& \rsect_{[S \cap S']} \shm.
\end{eqnarray}
Let $f\colon X \to Y$ be a morphism of smooth varieties, $Z\subset Y$ a
possibly singular closed subvariety, and set $S = f^{-1}(Z) \subset X$. 
Then there is an isomorphism
\begin{equation}
\label{eq:sectfSZ}
\doim f \rsect_{[S]} \shm \simeq \rsect_{[Z]} \doim f \shm.
\end{equation}
Let $Y$ be a closed smooth subvariety of $X$ of codimension $d$, and
denote by $j\colon Y \to X$ the embedding. 
Recall that Kashiwara's equivalence states that the functors $\shm\mapsto \dopb j \shm \shift{-d}$ and $\shn \mapsto \doim j \shn$ establish an equivalence between the category
$\Mod_\qcoh(\DY)$ of quasi-coherent $\DY$-modules, and the full abelian subcategory of
$\Mod_\qcoh(\DX)$ whose objects $\shm$ satisfy $\sect_{[Y]}\shm \isoto
\shm$. This extends to derived categories.
In particular, the functor
\begin{equation}
\label{eq:kasheq}
\doim j \colon \BDC_{\qcoh}(\DY) \to \BDC_{\qcoh}(\DX) \quad \text{is fully faithful,}
\end{equation}
and for $\shm\in\BDC_{\qcoh}(\DX)$ one has
\begin{equation}
\label{eq:sectY}
\rsect_{[Y]}\shm \simeq \doim j \dopb j \shm \shift{-d}.
\end{equation}

\subsection{Fourier-Laplace transform}\label{sse:fourier}

To $\varphi \in \sect(X;\OX)$ one associates the
$\DX$-module\footnote{Equivalently, $\DX e^{\varphi}$ is the sheaf
$\OX$ with the $\DX$-module structure given by the flat connection
$1 \mapsto d\varphi$.}
$$
\DX e^{\varphi} = \DX / \shi_{\varphi}, \qquad \shi_{\varphi}(V) = \{
P\in\DX(V) \colon P e^{\varphi} = 0 \},\quad \forall V\subset X\text{ open}.
$$
For $f\colon X \to Y$ a morphism of smooth algebraic
varieties, and $\psi \in \sect(Y;\OY)$, one has
\begin{equation}
\label{eq:expopb}
\dopb f \DY e^{\psi} \simeq \DX e^{\psi \circ f}.
\end{equation}
Let us denote by $\AffX$ the trivial line bundle on $X$, and by $t
\in \sect(\AffX;\O_{\AffX})$ its fiber coordinate.  Let $\piV \colon
\VV \to X$ be a vector bundle of finite rank, $\piA
\colon \AA \to X$ be the dual bundle, $\gamma_{\VV} \colon \VVAA \to \AffX$ be the natural pairing, and $\VV\from[p_1]\VVAA\to[p_2]\AA$ be the natural projections.
The Fourier-Laplace transform for $\D$-modules is the functor
\begin{align*}
   \dfourier[\VV] \colon \BDC_\qcoh(\DV) &\to  \BDC_\qcoh(\DA)\\
   \shn &\mapsto \doim{p_{2}{}} (\dopb{p_{1}} \shn \dtens 
   \dopb{\gamma_{\VV}} \shl_1).
\end{align*}
The Fourier-Laplace transform is involutive, in the sense that (cf~\cite[Lemma 7.1 and Appendix 7.5]{Kz-L})
\begin{equation}
\label{eq:fourierinvolutive}
\dfourier[\AA] \circ \dfourier[\VV] \simeq
\dopb{(-\id_{\VV})}.
\end{equation}
Let $f\colon \VV\to \WW$ be a morphism of vector bundles over $X$, and
denote by $\tr{f}\colon\MM\to\AA$ the transpose of $f$.  Then for any
$\shn\in\BDC_{\qcoh}(\DV)$ and $\shp\in\BDC_{\qcoh}(\DW)$ there are
natural isomorphisms\footnote{See the appendix for a proof.}
\begin{eqnarray}
\label{eq:oimfourier}
\dfourier[\WW]\doim f \shn &\simeq&
\dopb{\ttr{f}}\dfourier[\VV]\shn, \\
\label{eq:opbfourier}
\dfourier[\VV] \dopb{f} \shp &\simeq& \doim{\ttr{f}}
\dfourier[\WW] \shp.
\end{eqnarray}
If $X$ is viewed as a zero-dimensional vector bundle over itself, the
projection $\piV\colon\VV\to X$ and the zero-section $\iA \colon X\to
\AA$ are transpose to each other.  Hence \eqref{eq:opbfourier} 
gives for $\shm\in\BDC_{\qcoh}(\DX)$ and
$\shq\in\BDC_{\qcoh}(\DA)$ the isomorphisms\footnote{Note that isomorphism \eqref{eq:oimi} is the content of~\cite[Lemma
2.3]{D-M-S-S}, of which we have thus provided a more natural proof.}
\begin{eqnarray}
\label{eq:oimi}
\doim{\iA}\shm &\simeq& \dfourier[\VV]
\dopb\piV\shm, \\
\label{eq:opbi}
\dopb{\iA}\shq &\simeq&
\doim\piV\dfourier[\AA]\shq.
\end{eqnarray}

\section{Dwork cohomology}

Let $s\colon X\to\AA$ be a section of the vector bundle
$\piA\colon\AA\to X$ of rank $r$, and set $\tilde s =
\id_{\VV}\times_{X} s\colon \VV \to \VVAA$.  
Recall that $\gamma_{\VV}\colon \VVAA \to \AffX$ denotes the pairing, and let $F\in\sect(\VV;\OV)$ be
the function 
$$
F = t \circ \gamma_{\VV} \circ \tilde s.
$$
Let us denote by $S$ the reduced zero locus of
$s$, which is a possibly singular closed subvariety of $X$, and
by $j\colon S\to X$ the embedding. The geometric framework is thus
summarized in the commutative diagram with Cartesian squares
$$
\xymatrix{ X \ar[r]^{\iA} \ar@{}[dr]|-\square & \AA
\ar@{}[dr]|-\square & \VVAA \ar[l]_{{p_{2}}} \ar[r]^{\gamma_{\VV}} 
\ar[dr]^{p_1} & \AffX \\
\widetilde  S \ar[u]^{j_1} \ar[r]^{j_2} & X \ar[u]^{s} & \VV \ar[l]_{\piV}
\ar[u]^{\tilde s} & \VV \ar[l]_{\id_\VV}. }
$$
Then $j_1 = j_2 = j$ on $S = \widetilde S^{\text{red}}$.
Generalizing previous results of~\cite{Kz,M,A-S}, Theorem 0.2 of \cite{D-M-S-S} gives the following link between relative cohomology and Dwork cohomology

\begin{theorem}
For $\shm\in\BDC_{\qcoh}(\DX)$ there is an isomorphism
$$
\rsect_{[S]}\shm\shift{r} \simeq \doim\piV(\dopb\piV \shm \dtens \DV
e^{F}).
$$
\end{theorem}

Our aim here is to provide a more natural proof of this result.

\begin{proof}
For $\shm = \OX$, the statement reads
\begin{equation}
\label{eq:dwork}
\rsect_{[S]} \OX \shift{r} \simeq \doim\piV\DV e^{F}  .
\end{equation}
For a general $\shm\in\BDC_{\qcoh}(\DX)$, there are isomorphisms
$$
\rsect_{[S]}\shm \simeq \shm \dtens \rsect_{[S]}\OX \qquad \text{by }
\eqref{eq:sectM},
$$
and
$$
\doim\piV(\dopb\piV \shm \dtens \DV e^{F}) \simeq \shm\dtens \doim\piV\DV
e^{F} \qquad \text{by } \eqref{eq:projection}.
$$
It is thus sufficient to prove \eqref{eq:dwork}.
Setting $\shl = \dopb{\gamma_{\VV}} \DAffX e^{t}$, there is a chain of isomorphisms
\begin{align*}
\DV e^{F} & \simeq  \dopb{\tilde s}
\shl && \text{by }\eqref{eq:expopb} \\
& \simeq  \dopb{\tilde s} \shl \dtens \OV \\
& \simeq  \doim{p_{1}{}} \doim{\tilde s} (\dopb{\tilde s}
\shl \dtens \OV) && \text{by }\eqref{eq:fcircgoim}  \\
& \simeq  \doim{p_{1}{}} (\shl
\dtens \doim{\tilde s} \OV) && \text{by }\eqref{eq:projection} \\
& \simeq  \doim{p_{1}{}} (\shl
\dtens \doim{\tilde s} \dopb{\piV} \OX) \\
& \simeq  \doim{p_{1}{}} (\shl
\dtens \dopb{p_{2}} \doim s \OX) && \text{by }\eqref{eq:basechange} \\
& =  \dfourier[\AA] \doim s \OX.
\end{align*}
Hence we have
\begin{align*}
\doim{\piV{}}\DV e^{F} & \simeq  
\doim{\piV{}} \dfourier[\AA] \doim s \OX \\
& \simeq \dopb{\iA} \doim s \OX && \text{by }\eqref{eq:opbi},
\end{align*}
and to prove \eqref{eq:dwork} we are left to establish an isomorphism
\begin{equation}
\label{eq:jjrsect}
\dopb{\iA} \doim s \OX \simeq \rsect_{[S]}\OX\shift{r}.
\end{equation}
By \eqref{eq:kasheq}, this follows from the chain of isomorphisms
\begin{align*}
\doim{\iA} \dopb{\iA} \doim s \OX 
& \simeq \doim{\iA} \dopb{\iA} \doim s \dopb s \OA \\
& \simeq \rsect_{[\iA(X)]} \rsect_{[s(X)]} \OA \shift{2r} && \text{by
}\eqref{eq:sectY} \\
& \simeq \rsect_{[\iA(S)]} \rsect_{[\iA(X)]} \OA
\shift{2r} && \text{by }\eqref{eq:sectSS'} \\
& \simeq \rsect_{[\iA(S)]} \doim{\iA} \OX \shift{r} && \text{by
}\eqref{eq:sectY} \\
& \simeq \doim{\iA} \rsect_{[S]} \OX \shift{r} && \text{by
}\eqref{eq:sectfSZ}.
\end{align*}
\end{proof}

\begin{remark}
Kashiwara's equivalence allows one\footnote{This is done for example
in~\cite{B}.  For $S$ a singular closed subvariety of a smooth variety
$X$, the idea is to define $\Mod(\D_{S})$ as the full abelian
subcategory of $\Mod(\DX)$ whose objects $\shm$ satisfy
$\sect_{[S]}\shm\isoto\shm$.} to develop the theory of algebraic
$\D$-modules on possibly singular varieties, so that
the formulae stated in the previous section still hold.  In this framework,
\eqref{eq:jjrsect} is obtained by
\begin{align*}
\dopb{\iA} \doim s \OX & \simeq \doim j \dopb j \OX &&
\text{by }\eqref{eq:basechange} \\
& \simeq \rsect_{[S]}\OX\shift{r} && \text{by }\eqref{eq:sectY}.
\end{align*}
\end{remark}

\appendix
\section{Appendix}

\subsection{Base change and projection formulae}

The base change formula \eqref{eq:basechange} is proved in~\cite[Theorem~VI.8.4]{Bo} for $h$ a locally closed embedding\footnote{In the language of Gauss-Manin connections, the base change formula is stated in~\cite[\S~3.2.6]{An-Ba} for $h$ flat.}. Let us recall how to deal with the general case.

\begin{proof}[Proof of \eqref{eq:basechange}]
The Cartesian square \eqref{eq:cartsq}
splits into the two Cartesian squares
$$
\xymatrix@C=5em{ 
X' \ar[r]^{(f',h')} \ar[d]_{f' } & Y'\times X  \ar[r]^{p'_2} \ar[d]|{\id_{Y'}\times f}  \ar@{}[dl]|-\square & X \ar[d]_{f} \ar@{}[dl]|-\square \\ 
Y' \ar[r]^{(\id_{Y'},h)} & Y' \times Y \ar[r]^{p_2} & Y,
}
$$
where $p_2$ and $p'_2$ are the natural projections.
Since $(\id_{Y'},h)$ is a closed embedding, by~\cite{Bo} the base change formula holds for the Cartesian square on the left hand side. We are thus left to prove the base change formula for the Cartesian square on the right hand side. For $\shm\in\BDC_\qcoh(\DX)$, one has the chain of isomorphisms
\begin{align*}
\doim{(\id_{Y'}\times f)} \dopb{p_2'{}} \shm 
& \simeq \doim{(\id_{Y'}\times f)} ( \O_{Y'} \etens \shm ) \\
& \simeq  \O_{Y'} \etens \doim{f} \shm \\
& \simeq  \dopb{p_2} \doim f \shm,
\end{align*}
where $\etens$ denotes the exterior tensor product.
\end{proof}

Let us also recall, following~\cite{B}, how projection formula is deduced from base change formula.

\begin{proof}[Proof of \eqref{eq:projection}]
Consider the diagram with commutative square
$$
\xymatrix{
& X \ar[dl]_{\delta_X} \ar[r]^{f} \ar[d]^{\delta_f} & 
Y \ar@{}[dl]|-\square \ar[d]^{\delta_Y} \\
X\times X \ar[r]^{f''}
& X\times Y \ar[r]^{f'} & Y\times Y, 
}
$$
where $\delta_X$ and $\delta_Y$ are the diagonal embeddings, $\delta_f$ is the graph embedding, $f'= f\times \id_Y$, and $f'' = \id_X \times f$. Then there is a chain of isomorphisms
\begin{align*}
\doim f (\shm \dtens \dopb f \shn)
& \simeq \doim f \dopb{\delta_X}(\shm \etens \dopb f \shn) \\
& \simeq \doim f \dopb{\delta_X}\dopb{f''}(\shm \etens \shn) \\
& \simeq \dopb{\delta_Y}\doim{f'}(\shm \etens \shn)  && \text{by }\eqref{eq:basechange} \\
& \simeq \dopb{\delta_Y}(\doim f\shm \etens \shn) \\
& \simeq \doim f\shm \dtens \shn.
\end{align*}
\end{proof}

\subsection{Fourier-Laplace transform}

The formulae stated in section~\ref{sse:fourier} for the
Fourier-Laplace transform of algebraic $\D$-modules have their analogues
for the Fourier-Deligne transform of $\ell$-adic sheaves
(see~\cite{L} or \cite[\S III.13]{Kl-W}), and for the Fourier-Sato transform of
conic abelian sheaves (see~\cite{K-Ssom}). Apart from \cite{Kz-L}, we do not have specific references for the algebraic $\D$-module case. We thus provide  here some proofs.

\begin{proof}[Proof of \eqref{eq:oimfourier} and \eqref{eq:opbfourier}]
The following arguments are parallel to those in the proof of
\cite[Th\'eor\`eme 1.2.2.4]{L} or \cite[Proposition 3.7.14]{K-Ssom}.
Consider the diagram with Cartesian squares $$\xymatrix{ \AA
\ar@{}[dr]|-\square & \MM \ar[l]_{\tr{f}} & \\
\VVAA \ar[u]^{p_{2}} \ar[dr]_{p_{1}} & \VVMM \ar[l]_{\alpha}
\ar[u]^{r_{2}} \ar[r]^{\beta} \ar[d]_{r_{1}} & \WWMM \ar[ul]_{q_{2}}
\ar[d]_{q_{1}} \\
& \VV \ar[r]_{f} & \WW \ar@{}[ul]|-\square , }
$$
where the morphisms $p_i$, $q_i$, $r_i$, for $i=1,2$ are the natural projections.
Note that $\gamma_{\VV} \circ \alpha = \gamma_{\WW} \circ \beta$.
The isomorphism \eqref{eq:oimfourier} is obtained via the following chain of
isomorphisms\footnote{Note that these arguments still apply if one
replaces $\shl_1 = \DAffX e^{t}$ with an arbitrary quasi-coherent
$\DAffX$-module.  On the other hand, in order to prove
\eqref{eq:opbfourier} we will use the fact that the Fourier transform
is involutive.  }, where we set $\shl_1 = \DAffX e^{t}$.
\begin{align*}
\dopb{\ttr{f}}\dfourier[\VV]\shn & =
\dopb{\ttr{f}}\doim{p_{2}{}}(\dopb{p_{1}}\shn\dtens\dopb{\gamma_{\VV}}
\shl_1) \\
& \simeq
\doim{r_{2}{}}\dopb{\alpha}(\dopb{p_{1}}\shn\dtens\dopb{\gamma_{\VV}}
\shl_1) && \text{by }\eqref{eq:basechange} \\
& \simeq
\doim{r_{2}{}}(\dopb{\alpha}\dopb{p_{1}}\shn\dtens\dopb{\alpha}\dopb{\gamma_{\VV}}
\shl_1) && \text{by }\eqref{eq:fopbtens} \\
& \simeq
\doim{r_{2}{}}(\dopb{r_{1}}\shn\dtens\dopb{\alpha}\dopb{\gamma_{\VV}}
\shl_1) && \text{by }\eqref{eq:fcircgopb} \\
& \simeq
\doim{r_{2}{}}(\dopb{r_{1}}\shn\dtens\dopb{\beta}\dopb{\gamma_{\WW}}
\shl_1) && \text{by }\eqref{eq:fcircgopb} \\
& \simeq \doim{q_{2}{}} \doim{\beta}
(\dopb{r_{1}}\shn\dtens\dopb{\beta}\dopb{\gamma_{\WW}} \shl_1)
&& \text{by }\eqref{eq:fcircgoim} \\
& \simeq
\doim{q_{2}{}}(\doim{\beta}\dopb{r_{1}}\shn\dtens\dopb{\gamma_{\WW}}
\shl_1) && \text{by }\eqref{eq:projection} \\
& \simeq
\doim{q_{2}{}}(\dopb{q_{1}}\doim{f}\shn\dtens\dopb{\gamma_{\WW}} \shl_1) && \text{by }\eqref{eq:basechange} \\
& = \dfourier[\WW]\doim f \shn.
\end{align*}
Applying the functor $\dfourier[\MM]$ to the isomorphism \eqref{eq:oimfourier} with $\shn = \dfourier[\AA]\shq$,
$\shq\in\BDC_{\qcoh}(\DA)$, and using \eqref{eq:fourierinvolutive}, we get
$$
\doim f \dfourier[\AA] \shq \simeq
\dfourier[\MM] \dopb{\ttr{f}} \shq.
$$
The isomorphism \eqref{eq:opbfourier} is obtained from the one above
by interchanging the roles of $f$ and $\tr f$.
\end{proof}

\providecommand{\bysame}{\leavevmode\hbox to3em{\hrulefill}\thinspace}

\makelastpage


\begin{thebibliography}{10}


\bibitem{An-Ba} Y. Andr\'e and F. Baldassarri, \emph{De Rham cohomology of differential modules on algebraic varieties},  Progress in Mathematics, 189, Birkh\"auser, 2001.

\bibitem{A-S} A. Adolphson and S. Sperber, \emph{Dwork cohomology, de Rham
 cohomology, and hypergeometric functions},  Amer. J. Math.  \textbf{122}
 (2000),  no. 2, 319--348.

\bibitem{B} J. Bernstein, Lectures on algebraic $\D$-modules at
Berkeley, 2001 (unpublished).

\bibitem{Bo} A. Borel, \emph{Algebraic $\D$-modules}, Perspectives in
Mathematics, 2.  Academic Press, 1987.

\bibitem{D-M-S-S} A. Dimca, F. Maaref, C. Sabbah, M. Saito,
\emph{Dwork cohomology and algebraic $\D$-modules}, Math.  Ann. 
\textbf{318} (2000), no.  1, 107--125.

\bibitem{K0} M.~Kashiwara, \emph{Algebraic study of systems of partial
differential equations}, M\'em.  Soc.  Math.  France (N.S.) (1995),
no.~63, xiv+72, Kashiwara's Master's Thesis, Tokyo University 1970,
translated from the Japanese by A.~D'Agnolo and J.-P.~Schneiders.

\bibitem{K1} \bysame, \emph{D-modules and microlocal
calculus} (translated from the 2000 Japanese original by M. Saito),
Transl. of Math. Monographs \textbf{217}, A.M.S. (2003).

\bibitem{K-Ssom} M. Kashiwara and P. Schapira, \emph{Sheaves on
manifolds}, Grundlehren der Mathematischen Wissenschaften, 292,
Springer, 1990.

% \bibitem{K-Slaplace} \bysame, \emph{Integral transforms with
% exponential kernels and {L}aplace transform}, J. Amer.  Math.  Soc. 
% \textbf{10} (1997), 939--972.

\bibitem{Kz} N.M. Katz, \emph{On the differential equations satisfied by period
matrices},  Inst. Hautes \'Etudes Sci. Publ. Math. \textbf{35} (1968) 223--258. 

\bibitem{Kz-L} N.M. Katz and G. Laumon, \emph{Transformation de
Fourier et majoration de sommes exponentielles}, Inst.  Hautes
\'Etudes Sci.  Publ.  Math.  \textbf{62} (1985), 361--418.

\bibitem{Kl-W} R. Kiehl and R. Weissauer, \emph{Weil conjectures, preverse sheaves and $l$'adic Fourier transform}, Ergebnisse der Mathematik und ihrer Grenzgebiete, 42, Springer, 2001.

\bibitem{L} G. Laumon, \emph{Transformation de Fourier, constantes
  d'\'equations fonctionnelles et conjecture de Weil}, Inst. Hautes
\'Etudes Sci. Publ. Math. \textbf{65} (1987),  131--210. 

\bibitem{M} P. Monsky, \emph{Finiteness of de Rham cohomology},
  Amer. J. Math.  \textbf{94} (1972), 237--245.

\end{thebibliography}
\end{document}